\documentclass{amsart} 
\usepackage{amsfonts}
\usepackage{epsf}
\usepackage{graphicx}
\usepackage{amssymb}
\usepackage{amsmath}
\usepackage{amscd}
\input diagrams

\newcommand{\OnS}{{Ozsv\'ath and Szab\'o }}
\newcommand{\zee}{\mbox{$\mathbb Z$}}
\newcommand{\spinc}{\mbox{{spin}${}^c$ }} 
\newcommand{\Spinc}{\mbox{{\rm Spin}${}^c$}}
\newcommand{\s}{{\mathfrak s}} 
\newcommand{\sym}{{\mbox{Sym}}} 

\newcommand{\x}{{\mbox{\bf x}}}

\newcommand{\R}{{\mathfrak r}}
\newcommand{\im}{\mbox{Im}}
\newcommand{\coker}{\mbox{\rm coker}}
\newcommand{\Ext}{\mbox{\rm Ext}}
\newcommand{\longto}{\longrightarrow}

\newcommand{\cal}{\mathcal}
\newcommand{\mb}{\boldmath}

\newcommand{\balpha}{\mbox{\mb $\alpha$}}
\newcommand{\bbeta}{\mbox{\mb $\beta$}}
\newcommand{\sbar}{\underline{\s}}

\newtheorem{thm}{Theorem}
\newtheorem{cor}[thm]{Corollary}
\newtheorem{lemma}[thm]{Lemma}

\newtheorem{prop}[thm]{Proposition}
\newtheorem{rem}[thm]{Remark}

\begin{document}

\title{Heegaard Floer Homology of Mapping Tori II}

\author{Stanislav Jabuka, Thomas Mark}

\begin{abstract} We extend the techniques in a previous paper to
calculate the Heegaard Floer homology groups $HF^+(M, \s)$ for fibered
3-manifolds $M$ whose monodromy is a power of a Dehn twist about a
genus-1 separating circle on a surface of genus $g\geq 2$, where $\s$
is a non-torsion \spinc structure on $M$.
\end{abstract}
\maketitle

\section{Introduction}

After their introduction by Ozsv\'ath and Szab\'o in a remarkable
series of papers \cite{OS1,OS2,OS3,OS4}, Heegaard Floer homology
groups and related invariants of 3- and 4-dimensional manifolds have
rapidly become central tools in low-dimensional topology. This article
is concerned with the calculation of the Heegaard Floer groups
associated to certain fibered 3-manifolds, and can be seen as a
continuation of, or supplement to, the work in \cite{us}. In the
latter paper, the authors considered 3-manifolds $Y$ fibered over
$S^1$ with fiber a surface $\Sigma_g$ of genus $g>1$ and monodromy
given as certain combinations of Dehn twists along circles in
$\Sigma_g$. Specifically, let $\gamma$ and $\delta$ denote a pair of
nonseparating circles on $\Sigma_g$ that intersect transversely in a
single point, and let $\sigma$ denote a circle that separates
$\Sigma_g$ into components of genus $1$ and $g-1$. Letting $t_\gamma$,
$t_\delta$, and $t_\sigma$ denote the right-handed Dehn twists about
$\gamma$, $\delta$, and $\sigma$, the monodromies considered in
\cite{us} were those of the form $t_\gamma^nt_\delta^m$ for
$n,m\in\zee$ and $t^{\pm 1}_\sigma$. For technical reasons the
calculations of that paper did not apply to the case of other powers
of the separating twist; part of the purpose of this paper is to make
the necessary extensions in the arguments so as to include those
cases.

Seen in a broader context, the present work extends to more general
circumstances a result initially used in \cite{OSknot} and
reformulated and applied in \cite{us}. Namely, given a nullhomologous knot
$K$ in a closed oriented 3-manifold $Y$ we determine (under some
hypotheses) the Heegaard Floer homology groups $HF^+(Y_0, \s)$ of the
result $Y_0$ of 0-framed surgery along $K$ in a \spinc structure $\s$
whose first Chern class is not torsion, in terms of the ``knot
complex'' $CFK^{\infty}(Y, K)$. This result is stated formally in
Theorem \ref{0surgthm} below, and the calculation of the Heegaard Floer
groups for the fibered 3-manifolds mentioned above is given as an
application.

For an orientation-preserving diffeomorphism $\phi$ of a surface
$\Sigma_g$, we denote by $M(\phi)$ the mapping torus $\Sigma_g\times
[0,1]/(x, 1)\sim (\phi(x),0)$. For $\sigma$ a genus-1 separating
circle on $\Sigma_g$ as above, note that $M(t_\sigma^n)$ has the
homology of $\Sigma_g\times S^1$ for any $n\in\zee$. It follows that
the following two conditions uniquely determine a \spinc structure
$\s_k$ on $\Sigma_g$:
\begin{enumerate}
\item[(i)] $\langle c_1(\s_k), [c]\times S^1\rangle = 0$ for any class
$[c]\in H_1(\Sigma_g; \zee)$.
\item[(ii)] $\langle c_1(\s_k), [\Sigma_g]\rangle = 2k$
\end{enumerate}
According to the adjunction inequality for Heegaard Floer homology
\cite{OS2}, $HF^+(M(t_\sigma^n), \s) = 0$ unless $\s$ satisfies (i)
above, and also satisfies (ii) with $|k|\leq g-1$. We calculate here the
Heegaard Floer groups in all remaining cases except $k = 0$.

To state the result, let $X(g,d)$ denote the graded group whose
summand in degree $j$ is $H^{g-j}(\sym^d \Sigma_g; \zee)$. For a
group $G$, we let $G_{(j)}$ denote the graded group isomorphic to $G$
with grading concentrated in degree $j$, and if $H$ is a graded group
let $H[j]$ denote $H$ with the grading shifted by $j$.

The following is proved in sections \ref{proofsec} and
\ref{lefttwist}, and incorporates Theorem 1.3 of \cite{us} as the cases $n
= \pm 1$.

\begin{thm}\label{mainthm} For $n\neq 0$, let $M(t_\sigma^n)$ denote the
mapping torus of the $n$-th power of the right-handed Dehn twist
around a genus 1 separating curve on a surface $\Sigma_{g}$ ($g\geq
2$). Let $\s$ be a nontorsion \spinc structure on $M(t_\sigma^n)$.
Then the Heegaard Floer homology $HF^+(M(t_\sigma^n), \s)$ is trivial
unless $\s$ satisfies (i) and (ii) above. If these conditions hold we
have an isomorphism of relatively graded groups
\begin{eqnarray}
HF^+(M(t_\sigma^n), \s) &=& X(g-1,d-1)\otimes H^*(S^1\sqcup
S^1)[\varepsilon(n)]\oplus
\Lambda^{2g-2-d}_{(g-d)}H^1(\Sigma_{g-1})\nonumber\\
&& \oplus \bigoplus_{p = 1}^{d} \left[ \Lambda^{2g-2-d+p}_{(g-d-p+1)}
H^1(\Sigma_{g-1})[\varepsilon(n)] \otimes H_*(\coprod_{|n|-1}
S^{2p-1})\right],\label{answer}
\end{eqnarray}
where $d = g-1-|k|$, and where $\varepsilon(n) = 0$ if $n>0$ and
$\varepsilon(n) = -1$ if $n<0$.
\end{thm}

Note that while the Heegaard Floer groups admit further algebraic
structure, in particular an action by the polynomial ring $\zee[U]$,
the methods in this paper do not give information about that action:
the isomorphism above is of $\zee$-modules only.

The following should be compared with results of Seidel \cite{seidel} and
Eftekhary \cite{eaman} on the symplectic Floer homology of surface
diffeomorphisms.

\begin{cor}
There is an isomorphism of relatively graded groups
\[
HF^+(M(t_\sigma^n), \s_{g-2}) = \left\{ \begin{array}{ll} H^*(\Sigma_g,
C) & n>0 \\ H^*(\Sigma_g \setminus C) & n<0\end{array}\right.
\]
where $C$ denotes the union of $|n|$ pairwise disjoint pushoffs of the
separating circle $\sigma\subset\Sigma_g$, and the right-hand side
above denotes singular cohomology with coefficients in $\zee$.
\end{cor}

\begin{proof} We consider the case $n>0$. In the statement we take $k
= g-2$, so that $d = 1$. Since $X(g-1,0) = \zee_{(g-1)}$, the first
line of \eqref{answer} becomes
$\zee_{(g-1)}\otimes(\zee_{(0)}^2\oplus\zee_{(1)}^2)\oplus
\Lambda^{2g-3}H^1(\Sigma_{g-1})_{(g-1)} =
\zee_{(g-1)}^2\oplus\zee_{(g)}^2\oplus \zee^{2g-2}_{(g-1)}$. The
second line of \eqref{answer} collapses to
$\zee_{(g-1)}\otimes(\zee_{(0)}^{n-1}\oplus \zee^{n-1}_{(1)}) =
\zee_{(g-1)}^{n-1}\oplus\zee_{(g)}^{n-1}$. Thus
\[
HF^+(M(t_\sigma^n), \s_{g-2}) \cong \zee^{n+1}_{(g)} \oplus
\zee^{2g + n -1}_{(g-1)}.
\]
On the other hand, since $\Sigma_g\setminus C$ is the disjoint union
of a punctured torus, a punctured surface of genus $g-1$, and $n-1$
annuli, we have
\[
H^*(\Sigma_g, C) \cong \zee^{n+1}_{(2)}\oplus \zee^{2g+n -1}_{(1)}
\]
which agrees with the above modulo a shift in grading. The case $n<0$
follows similarly.
\end{proof}

In the next section we state and prove Theorem \ref{0surgthm} on the
Floer homology of a manifold obtained by 0-framed surgery on a
nullhomologous knot; in Section \ref{calcsec} we prove Theorem
\ref{mainthm}.

\section{Floer Homology of 0-Surgeries}

The method we will use to make our calculation is an adaptation to
slightly more general circumstances of the one used by \OnS in
\cite{OSknot} to determine the Heegaard Floer homology of the mapping
torus of a single Dehn twist about a nonseparating curve on a surface of
genus $g\geq 2$, in nontorsion \spinc structures.

Let $K\subset Y$ be a nullhomologous knot in a closed oriented 3-manifold
$Y$, and let $Y_n(K)$ denote the result of $n$-framed surgery along
$K$. Fix a \spinc structure $\s$ on $Y$, and for integers $n>0$ and $k\in
\{0,\ldots,n-1\}$ let $HF^+(Y_n(K), k)$ denote the Heegaard Floer
homology of the surgered manifold in the \spinc structure $\s_k$
defined as follows. Integer surgery corresponds to a cobordism $W_n$
from $Y$ to $Y_n(K)$, comprising a single 2-handle addition. Fixing a
Seifert surface $F$ for $K$ we obtain a closed surface $\hat{F}$ in
$W_n$ by capping off $F$ using the core of the 2-handle . Now let $\s_k\in \Spinc(Y_n(K))$ be the \spinc structure
cobordant to $\s$ by a \spinc structure $\R$ on $W_n$ having
\[
\langle c_1(\R), [\hat{F}]\rangle = n - 2k.
\]

Note that if, as will always be the case here, the \spinc structure on
$Y$ is torsion and $H_2(Y;\zee)$ is torsion-free, then $\s_k$ is
independent of the choice of $F$. 

According to \cite{OS2}, there is an integer surgeries long exact
sequence
\begin{equation}
\cdots\to HF^+(Y_0(K),[k])\to HF^+(Y_n(K),k)\to HF^+(Y)\to
\cdots
\label{les}
\end{equation}
Here the first term denotes the sum of Floer homology groups
over \spinc structures in the fiber of $\s_k$ under a certain
surjective map $\Spinc(Y_0(K))\to \Spinc(Y_n(K))$ (see \cite{OS2}).
Our object is to use knot Floer homology to understand each term in
this sequence, and therefore we quickly review the relevant facts
about knot Floer homology (for details, see \cite{OSknot},
\cite{OSknot2}).

Given $K\subset Y$, let $E$ denote the torus boundary of a regular
neighborhood of $K$. One can then find a Heegaard surface for $Y$ of
the form $E\# \Sigma_{g-1}$, with attaching circles $\balpha =
\alpha_1,\ldots,\alpha_g$ and $\bbeta = \beta_1,\ldots,\beta_g$ where
$\beta_1\subset E$ is a meridian for $K$ and $(\Sigma,
\balpha, \beta_2,\ldots,\beta_g)$ is a Heegaard diagram for the knot
complement $Y\setminus K$. Let $w$
and $z$ denote a pair of basepoints, one on each side of the meridian
$\beta_1$. The data $(E\#\Sigma, \balpha,\bbeta, w)$ together with a choice
of \spinc structure $\s$ on $Y$ can be used to define the Heegaard
Floer chain groups $CF^\infty(Y, \s)$. The additional basepoint $z$,
along with a choice of ``relative \spinc structure''
$\underline{\s}\in Spin^c(Y_0(K))$ lifting $\s$ gives rise to a
filtration $\cal F$ on $CF^\infty(Y,\s)$. The ``knot chain complex''
$CFK^\infty(Y,K,\sbar)$ is this filtered complex.

More concretely, we fix a Seifert surface $F$ for $K$: then $F$
specifies the zero-framing on $K$, and can be capped off to a closed
surface $\hat{F}$ in the zero-surgery $Y_0(K)$. The generators of
$CFK^\infty(Y,K,F)$ are triples $[\x,i,j]$, where $\x$ denotes an
intersection point between the $g$-dimensional tori $T_\alpha =
\alpha_1\times\cdots\times \alpha_g$ and $T_\beta = \beta_1\times
\cdots\times\beta_g$ in the symmetric power $\sym^g(E\#\Sigma)$, and
$i$ and $j$ are integers. The point $\x$ along with the basepoint $w$
determine a \spinc structure $\s_w(\x)$ on $Y$ as well as a relative
\spinc structure $\underline{\s}_w(\x)$; we require that ${\s}_w(\x) =
{\s}$. Furthermore, $i$ and $j$ are required to satisfy the equation
\[
\langle c_1(\sbar_w(\x)), [\hat{F}]\rangle = 2(j-i).
\]
 In this notation, the filtration
$\cal F$ is simply ${\cal F}([\x,i,j]) = j$; changing the Seifert
surface $F$ shifts $\cal F$ by a constant.

The boundary map $\partial^\infty$ in $CFK^\infty$ is
defined by counting holomorphic disks in $\sym^g(E\#\Sigma)$ and can
only decrease the integers $i$ and $j$. Thus, for example, the subgroup 
$C\{i<0\}$ of $CFK^\infty$ generated by those $[\x,i,j]$ having $i<0$
is a subcomplex, and indeed is simply $CF^-(Y,\s)$ with an additional
filtration. The quotient of $CFK^\infty$ by $C\{i<0\}$ is written
$C\{i\geq 0\}$, and is a filtered version of $CF^+(Y,\s)$. We will use
other similar notations to indicate other sub- or quotient complexes
of $CFK^\infty$. In particular, $CFK^{0,*}$ is by definition the quotient
complex $C\{i = 0\}$, and $\widehat{HFK}$ is the homology of the
graded object associated to the filtration ${\cal F}$ of
$CFK^{0,*}$. We denote by $\widehat{HFK}(Y,K;j)$ the summand of
this group supported in filtration level $j$ (typically suppressing
the \spinc structure from the notation).

As an additional piece of structure, we have a natural chain endomorphism
$U$ on $CFK^\infty$ given by $U:[\x,i,j]\mapsto [\x,i-1,j-1]$.

\OnS prove (Theorem 4.4 of \cite{OSknot}):

\begin{thm} \label{hello} For all
sufficiently large positive $n$, there exists a $U$-equivariant
isomorphism of chain complexes
\[
{^b\Psi^+} : CF^+ (Y_n , k) \rightarrow C\{ i \ge 0 \mbox{ or } j \ge k \} 
\]
\end{thm}

In particular, $HF^+(Y_n(K), k)$ is given by the homology of the
portion of the knot complex indicated on the right-hand side. It is
important to note that the proof of the above theorem shows that the
stated identification is induced by a chain map coming from a
particular \spinc structure on the cobordism $-W_n$ connecting
$Y_n(K)$ to $Y$, namely the structure $\R_0\in\Spinc(-W_n)$ satisfying
\begin{equation}\label{specialstr}
\langle c_1(\R_0),[\hat{F}]\rangle = n
- 2k.
\end{equation}
By contrast, the homomorphism $HF^+(Y_n(K),k)\to HF^+(Y)$
in the long exact sequence (\ref{les}) is given by the sum of the
homomorphisms induced by all \spinc structures on $-W_n$ that restrict
to the given \spinc structures on $Y_n(K)$ and $Y$. We will return to
this point shortly.

Examining the sequence (\ref{les}) again, we can now understand two of
the groups appearing in terms of the knot chain complex: the theorem
above identifies $HF^+(Y_n,k)$, while $HF^+(Y)$ is simply
$H_*(C\{i\geq 0\})$, the homology of the complex obtained by
forgetting the filtration. Note that under these
identifications, the map $HF^+(Y_n,k)\to HF^+(Y)$ induced by $\R_0$ as
above corresponds to the map on homology induced by the natural
projection $C\{i\geq 0\mbox{ or } j\geq k\} \to C\{i\geq 0\}$.
Supposing the connecting homomorphism $HF^+(Y)\to
HF^+(Y_0(K), [k])$ to be trivial, and ignoring the issue raised in the
previous paragraph, we anticipate an isomorphism
\[
HF^+(Y_0(K), [k]) \cong H_*(C\{i<0\mbox{ and } j\geq k\}).
\]
Indeed, it is implicit in \cite{OSknot} and described concretely in
\cite{us} that if $HF^+_{red}(Y) = 0$ (and certain other, less
important, hypotheses hold) then the connecting homomorphism is in
fact trivial. We will show that the isomorphism holds even in certain
cases when the connecting homomorphism is nontrivial.

\begin{thm} \label{0surgthm} Fix a torsion \spinc structure $\s$ on $Y$,
and let $HF^+(Y, \s; A)$, $HF^+(Y_n(K), k; A)$, and $HF^+(Y_0(K), [k];
A)$ denote Floer homology groups in \spinc structures as above, with
coefficients in a ring $A$. Assume that $k$ is nonzero, and make the
following additional assumptions:
\begin{enumerate}
\item[1.] $HF^+(Y, \s, A)$ is a free $A$-module.
\item[2.] The filtration on
$HF_{red}(Y)$ is proportional to the degree: specifically, there is a
constant $c$ such that the part of $HF_{red}(Y)$ that lies in
filtration level $j$ is supported in absolute grading $j + c$. 
\item[3.] For $F: HF^+(Y_n(K),k;A) \to HF^+(Y;A)$ the homomorphism in
the long exact sequence \eqref{les}, we have
\[
\Ext^1_A(\ker(F), \coker(F)) = 0.
\]
\item[4.] The reduced homology $HF_{red}(Y_n(K), k,A)$ is supported in
degrees at most $k + c - d$, where $d$ is the degree shift induced by the
\spinc structure $\R_0$ as in \eqref{specialstr}.
\end{enumerate}
Then there is an identification of $A$-modules
\begin{equation}
\label{homologyident}
HF^+(Y_0(K), [k]; A) \cong H_*(C\{i<0\mbox{ and } j\geq k\}; A).
\end{equation}

\end{thm}

In particular hypotheses 1 and 3 of the theorem hold if $A$ is a field. If $A = 
\zee$, these hypotheses hold if, for example, $HF^+(Y)$ and $\coker(F)$ are 
torsion-free groups.

For the proof we need to analyze the surgery exact sequence:
\begin{equation}\label{surgseq}
\cdots\rTo^G HF^+(Y_0(K),[k])\rTo^H HF^+(Y_n,k)\rTo^F HF^+(Y)\rTo^G
\cdots
\end{equation}
In particular, as mentioned above, we study the homomorphism $F$.
Recall that since the \spinc structures $\s$ and $\s_k$ are torsion,
the corresponding Heegaard Floer homology groups admit a
(rational-valued) grading that lifts the natural relative $\zee$
grading. A \spinc structure $\R$ on $-W_n$ extending $\s$ restricts to
$\s_k$ on $Y_n(K)$ if and only if it satisfies
\[
\langle c_1(\R), [\hat{F}]\rangle = n - 2k + 2xn
\]
for some integer $x\in\zee$, where $\hat{F}$ denotes the capped-off
Seifert surface in $-W_n$ as before. Recall that we assume
$k\in\{1,\ldots, n-1\}$. Now, $F$ is the sum of the maps $F_{-W_n, \R}$
for $\R$ corresponding to all values of $x$. According to \cite{OS3},
the map $F_{\R}: HF^+(Y_n(K), k)\to HF^+(Y)$ induced by $\R$ shifts
degree by the quantity
\[
\deg(F_{\R}) = -nx^2 - (n-2k)x + \frac{n-(n-2k)^2}{4n}.
\]
The maximum value of this degree occurs when $x$ is the closest
integer to $-\frac{1}{2} + \frac{k}{n}$, which, given our assumption
on $k$, is $x = 0$. Comparing with \eqref{specialstr} we see that the
\spinc structure $\R_0$ corresponding to $x = 0$ is both the ``leading
order term'' (i.e., the homogeneous part with maximal degree) in the
homomorphism $F$ and also the \spinc structure inducing the
identification of $HF^+(Y_n(K),k)$ with $H_*(C\{i\geq 0\mbox{ or }
j\geq k\})$. 

With the above in mind, we write (as in \cite{OSknot}) $F = f_1 + f_2$
where $f_1$ is the homogeneous part with highest degree (that is, $f_1
= F_{-W_n, \R_0}$) and $f_2$ is the sum of all lower-degree parts of
$F$. Note that the highest-degree part of $f_2$
has degree equal to $\deg(f_1) - \min\{2k, 2(n-k)\}$.

\begin{lemma}\label{rinvlemma} Suppose there exists a homomorphism $R:
\im(f_2)\to HF^+(Y_n(K),k)$ satisfying $f_1\circ R = id$ on $\im(f_2)$
(that is, $R$ is a partially-defined right inverse for $f_1$). Then
there is an isomorphism $g: HF^+(Y_n(K),k)\to HF^+(Y_n(K),k)$ such
that
\[
f_1 = F\circ g.
\]
In particular, $\coker(F) = \coker(f_1)$, and $\ker(F) \cong
\ker(f_1)$ (via the map $g$).
\end{lemma}

\begin{proof}
Since $F = f_1 + f_2$, we have
\[
F = f_1( 1 + Rf_2).
\]
From the remarks on the degree shift above, the composition $Rf_2$ is
strictly decreasing in degree. Since $HF^+(Y_n(K),k)$ is trivial in
sufficiently low degrees, the sum
\[
g = (1 + Rf_2)^{-1} = \sum_{n\geq 0} (-Rf_2)^n
\]
is finite, and hence $g$ is the desired isomorphism.
\end{proof}

If we can construct the right inverse $R$ as in the lemma, Theorem
\ref{0surgthm} will follow easily. Indeed, let us identify
$HF^+(Y_n(K),k)$ with $H_*(C\{i\geq 0 \mbox{ or } j\geq k\})$ as in
Theorem \ref{hello}: then we've seen that $f_1$ is identified with the
natural projection

\[
\pi_*: H_*(C\{i\geq 0 \mbox{ or } j\geq k\})\to H_*(C\{i\geq 0\}) = HF^+(Y).
\]
Hence we have a commutative diagram{\small
\[
\begin{diagram}
\cdots & \rTo & HF^+(Y_0(K)) & \rTo & HF^+(Y_n(K)) & \rTo^F & HF^+(Y) &
\rTo & \cdots \\
      &       &  \dDashto    &      & \dTo >1+ Rf_2 &     & \dTo>id  &&\\
\cdots & \rTo & H_*(C\{i<0 \mbox{ and } j\geq k\}) & \rTo & 
H_*(C\{i\geq 0 \mbox{ or } j\geq k\}) & \rTo^{\pi_*} & H_*(C\{i\geq 
0\}) & \rTo & \cdots
\end{diagram}
\]}
where the solid vertical arrows are isomorphisms and the top row is 
the surgery long exact sequence \eqref{surgseq}. We have not 
constructed a map corresponding to the dashed arrow, but the above 
diagram allows us to write another:
\begin{equation}\label{commdiag}
\begin{diagram}
0 & \rTo & \coker(F) &\rTo& HF^+(Y_0(K)) &\rTo & \ker(F)& \rTo &0\\
&        &  \dTo>{id}     &    &              &     &  \dTo<{1 + Rf_2}  & \\
0 & \rTo & \coker(\pi_*) & \rTo & H_*(C\{i<0 \mbox{ and } j\geq k\}) &
\rTo & \ker(\pi_*) & \rTo & 0
\end{diagram}
\end{equation}
Again the vertical arrows are isomorphisms, so that
$HF^+(Y_0(K))$ and $ H_*(C\{i<0 \mbox{ and } j\geq k\})$ are both
extensions of $\coker(F)$ by $\ker(F)$. Under hypothesis 3 of the
theorem such an extension is unique, which completes the proof of
Theorem \ref{0surgthm}.

It remains to construct the right inverse $R$.

%
%
%

\begin{lemma}\label{rinvlemma2} Under hypotheses (1), (2), and (4) of
Theorem \ref{0surgthm} there exists a right inverse $R$ for $f_1$ that
is defined on $\im (f_2)$.
\end{lemma}

Note that the only place that hypothesis 3 of the theorem is required
is at the last step of its proof (see above).

\begin{proof} Since $\im(f_1)\subset HF^+(Y)$ is free by hypothesis, we can
find a right inverse $R: \im(f_1)\to HF^+(Y_n(K), k)$. We need to
check that $R$ is defined on $\im(f_2)$, i.e., that $\im(f_2)\subset
\im(f_1)$. We identify $HF^+(Y)$ with $H_*(C\{i\geq 0\})$ and
$HF^+(Y_n(K), k)$ with $H_*(C\{i\geq 0\mbox{ or } j\geq k\})$, so that
$f_1$ corresponds to the projection $\pi_*$ as above.

Recall that because of 
the structure of the chain complex $CF^{\infty}(Y)$ the image of 
the action of $U^r$ on $HF^+(Y)$ is independent of $r$ for 
sufficiently large $r$. By definition the reduced Floer homology group
$HF_{red}(Y)$ is the quotient of $HF^+(Y)$ by $\im(U^r)$ for any such
large $r$: in other words there is an exact sequence
\[
0 \rTo \im(U^r) \rTo HF^+(Y) \rTo HF_{red}(Y)\rTo 0.
\]
The knot filtration induces a filtration on $HF_{red}$.

Now, since the 
chain complexes $C\{i\geq 0\}$ and $C\{i\geq 0\mbox{ or } j\geq k\}$
are identical for sufficiently large degrees, it is clear that $f_1 =
\pi_*$ maps onto $\im(U^r)$ for $r>>0$. We claim that $f_1$ also maps
onto $H_{red}^{<k}(C\{i\geq 0\})$, where $H_{red}^{<k}$ denotes that
portion of the reduced homology that lies in filtration level $j=k-1$
or below. To see this, it suffices to show that the connecting
homomorphism $\delta: H_*(C\{i\geq 0\})\to H_*(C\{i<0\mbox{ and }j\geq
k\})$ is trivial on the indicated group. But for $x\in C\{i\geq 0\}$ a
cycle, $\delta x$ is given by the portion of $\partial^{\infty} x$
that lies in $C\{i< 0\mbox{ and } j\geq k\}$, where $\partial^\infty$
is the boundary map in $CFK^\infty$. The statement follows since
$\partial^\infty$ is nonincreasing in $j$.

 From the
discussion above, $f_1$ maps onto all summands of $HF^+(Y)$ that lie
in degree less than $k + c$ as well as $\im(U^r)$ (in all degrees). By
$U$-equivariance $f_2$ maps $\im(U^r)$ into $\im(U^r)$, so we need
only check that $f_2$ maps $HF_{red}(Y_n(K),k)$ into $\im(f_1)$. But
this follows immediately from the facts that $f_2$ has degree strictly
less than that of $f_1$ and that $HF_{red}(Y_n)$ is supported in degrees
$\leq k + c - \deg(f_1)$.
\end{proof}

\section{Calculation for Separating Twists}\label{calcsec}

\subsection{Preliminaries: Calculation of Knot Homology}

\begin{figure}[htb!] 
\centering
\includegraphics[width=12cm]{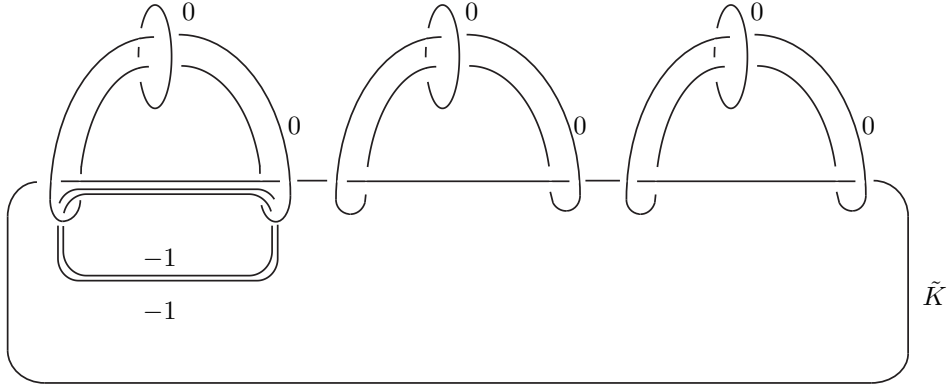}
\put(5,30){$\tilde{K}$}
\put(-60,138){$0$}
\put(-168,138){$0$}
\put(-275,138){$0$}
\put(-18,95){$0$}
\put(-127,95){$0$}
\put(-235,95){$0$}
\put(-290,25){$-1$}
\put(-290,45){$-1$}
\caption{The knot $\tilde K$ in the case $g=3$ and $n=2$.}  \label{sepknot2}
\end{figure}

A surgery diagram for the mapping torus $M(t_\sigma^n)$ is obtained
from one for $\Sigma_g\times S^1$ by adding $n$ parallel copies of the
separating curve $\sigma$ with surgery coefficient $-1$ (here and 
subsequently we assume $n>0$; the case of negative $n$ is entirely 
parallel and will be described later).
In particular, $M(t_\sigma^n)$ can be obtained from 0-surgery along
a knot $\tilde K$ in a connected sum $\#^{2g-2}(S^1\times S^2) \# M_n$,
where $M_n$ is described below. The situation is illustrated for $g = 3$
and $n = 1$ in Figure \ref{sepknot2}.

\begin{figure}[htb!] 
\centering
\includegraphics[width=4cm]{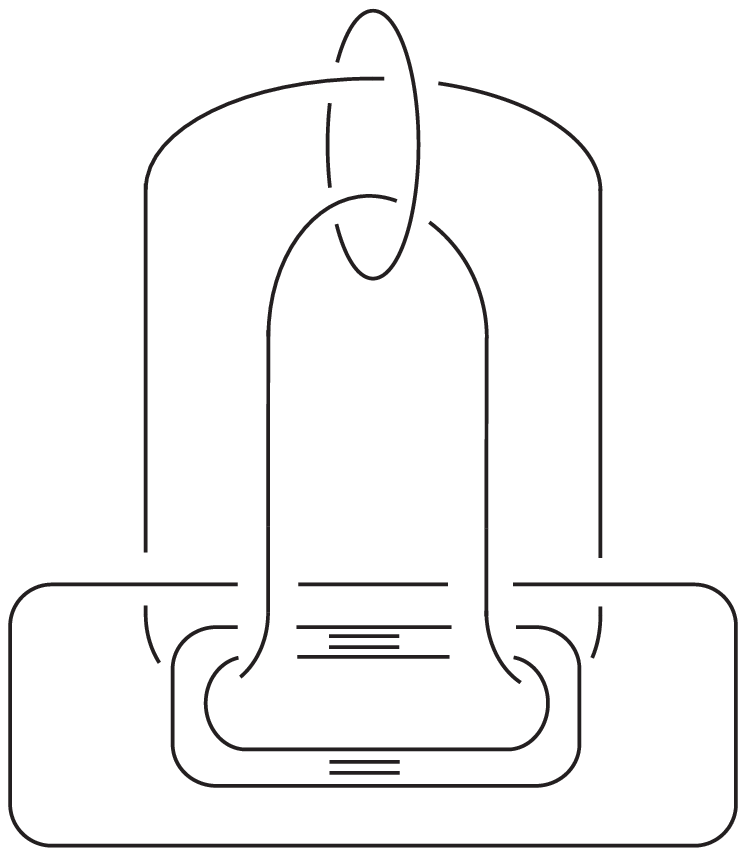} \quad \quad \quad \quad \quad \quad  \quad \quad 
\includegraphics[width=4cm]{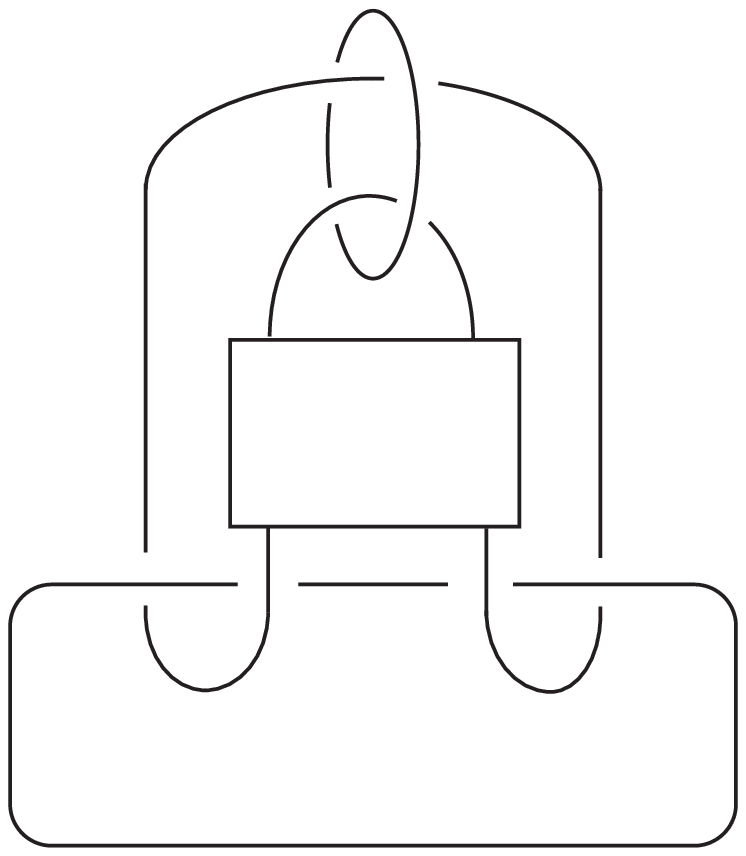}
\put(-62,62){$2n$}
\put(-50,127){0}
\put(-266,127){0}
\put(-19,80){0}
\put(-231,80){0}
\put(3,18){K}
\put(-209,18){K}
\put(-275,18){-1}
\put(-262,-25){(a)}
\put(-62,-25){(b)}
\caption{Two equivalent pictures for the knot $K$ in $M_n$. In
picture (a) there are $n$ parallel $-1$ circles. The $2n$ in
picture (b) indicates $2n$ positive half-twists. } \label{sepknot}
\end{figure}

We can think of $\tilde K$ as a connected sum of knots: write $B(p,q)$
for the knot given as the third component of the Borromean rings after
performing surgery on the other two components with surgery
coefficients $p$ and $q$. Then $\tilde K = \#^{g-1}B(0,0)\# K$, where
$K$ is the knot indicated in Figure \ref{sepknot}.

We denote by $M_n$ the manifold containing $K$, so $M_n$ is obtained
by performing 0-framed surgery on two components of the Borromean
rings and $-1$ framed surgery on each of $n$ parallel copies of the
third component. By 
blowing down all the $-1$ circles, we can see $M_n$ as the result of 
$0$-surgery on both components of an \lq\lq $n$-clasped Whitehead link\rq\rq 
as shown in figure \ref{sepknot}(b) (the $2n$ in the box denotes $2n$ 
positive half-twists). As a preliminary to the 
calculation of the knot homology of $K$, we find the homology 
$\widehat{HF}(M_n)$.

Let $Z$ denote the 3-manifold that results from changing the surgery
coefficient from $0$ to $-1$ on the smaller $0$-framed circle in the
diagram for $M_n$. An easy isotopy shows that after blowing down this
$-1$ circle, $Z$ is given as 0-framed surgery on the pretzel knot
$P(-2n + 1, -1, -1)$ (here we follow the notation conventions of
\cite{OS4}, whereby $P(1,1,1)$ is the right-handed trefoil). It
follows from the calculations in section 8 of \cite{OS4} that for
the torsion \spinc structure $\s_0$,
\[
HF^+_k(Z; \s_0) = \left\{\begin{array}{ll} \zee & k \equiv 
1/2 \mbox{ mod $\zee$ and $k\geq 3/2$}\\ \zee^n & k = 1/2 \\ 0 & 
\mbox{else}\end{array}\right.
\]
and that the image of $U$ in degree $1/2$ is nontrivial.

\begin{rem} For the rest of this section and throughout the next,
we take coefficients in the ring $A = \zee$. In particular it suffices
for the verification of hypothesis (3) of Theorem \ref{0surgthm} to
show that $\coker(F)$ is a free abelian group. 
\end{rem}

With this 
information, the long exact sequence in Floer homology connecting 
$S^1\times S^2$, $Z$, and $M_n$ immediately gives
\[
HF^+_k(M_n; \s_0) = \left\{\begin{array}{ll} \zee^2 & k \equiv 1 \mbox{
mod $\zee$ and $k\geq 1$}\\ \zee^{n+1} & k = 0 \\ 0 &
\mbox{else}\end{array}\right.
\]
From this we infer $\widehat{HF}(M_n) = \zee^{n+1}_{(1)}\oplus 
\zee^{n+1}_{(0)}$.

We are interested in the knot Floer homology $\widehat{HFK}(M_n, K)$. The
following was proved in \cite{us}:

\begin{lemma}\label{M1Klemma} The knot Floer homology groups for $(M_1, K)$ are given 
by
\[
\widehat{HFK}(M_1,K; j) = \left\{\begin{array}{ll} \zee_{(1)} & j = 
1\\ \zee^3_{(0)}\oplus \zee_{(1)} & j = 0 \\ \zee_{(-1)} & j = -1 \\ 
0 & \mbox{otherwise}\end{array}\right.
\]
The spectral sequence that calculates $\widehat{HF}(M_1)$ from 
$\widehat{HFK}(M_1,K)$ collapses at the $E_2$ level, and the only 
nontrivial differential is a surjection $d_1: \zee^3_{(0)}\to\zee_{(-1)}$.
\end{lemma}

We now calculate the result for the case of general $n$:

\begin{prop}\label{MnKprop} The knot Floer homology groups of $K$ in $M_n$ are given 
by
\[
\widehat{HFK}(M_n, K; j) = \left\{ \begin{array}{ll} \zee_{(1)} & j = 
1\\ \zee^{n+2}_{(0)} \oplus \zee^n_{(1)} & j = 0 \\ \zee_{(-1)} & 
j = -1 \\ 0 & \mbox{otherwise} \end{array}\right.
\]
The only nontrivial differential in the spectral sequence converging 
to $\widehat{HF}(M_n)$ is a surjection $d_1: \zee^{n+2}_{(0)}\to 
\zee_{(-1)}$.
\end{prop}

\begin{proof} 
We proceed inductively: suppose $n\geq 2$. We look at a surgery
sequence arising from figure \ref{sepknot}(a). Choose one of the $-1$
circles in the picture for $M_n$ and let $A$ denote the 3-manifold
obtained by changing the $-1$ to 0. Then the surgery sequence appears
as
\[
\cdots \to \widehat{HFK}(M_{n-1}, K; j)\to \widehat{HFK}(M_n, K ; j) 
\to \widehat{HFK}(A, K; j)\to \cdots
\]
Now $(A,K)$ is unknotted since $K$ can slide over the 0-framed 
circle. Therefore the filtration induced by $K$ on 
$\widehat{CF}(A)$ is trivial, and hence $\widehat{HFK}(A, K)$ is 
supported in level $j = 0$. It follows immediately that the group 
$\widehat{HFK}(M_n, K; j)$ has the claimed form when $j\neq 
0$.

For the case $j = 0$, note that the calculation so far together with 
the structure of $\widehat{HF}(M_n)$ already imply that 
$\widehat{HFK}(M_n, K; 0)$ is supported in degrees 0 and 1, by 
consideration of the spectral sequence for $\widehat{HF}$. It then 
follows that the only nontrivial differential in the latter spectral 
sequence is $d_1$, and in fact $d_1: \widehat{HFK}(M_n, K; 0)\to 
\widehat{HFK}(M_n, K, -1) \cong \zee$ must be surjective. An argument 
similar to the case $n = 1$ (see \cite{us}) together with our inductive
knowledge of $\widehat{HFK}(M_{n-1},K)$ shows that $d_1:
\widehat{HFK}(M_n, K, 1)\to \widehat{HFK}(M_n, K; 0)$ is trivial. The
proposition follows from this and the fact that $\widehat{HF}(M_n) =
\zee^{n+1}_{(1)} \oplus \zee^{n+1}_{(0)}$.
\end{proof}

It will be convenient in what follows to write $\widehat{HFK}(M_n, 
K)$ as
\[
\widehat{HFK}(M_n,K) \cong \Lambda^*H^1(\Sigma_1)\oplus H^*(\coprod_n S^1).
\]
In the above, the grading on $\Lambda^*H^1(\Sigma_1)$ is
``centered,'' meaning that $\Lambda^iH^1(\Sigma_g)$ is considered to
have grading $i-g$. The grading on the second factor above is the
usual homological grading. The filtration is not evident from the
notation, however we see that it is equal to the (centered) grading on
the first factor while the second factor lies in filtration level 0.

In this notation, we can express the single nontrivial differential 
in the spectral sequence for $\widehat{HF}(M_n)$ as the map 
$\Lambda^1H^1(\Sigma_1)\to \Lambda^0H^1(\Sigma_1)$ given by contraction with 
a generator $\gamma$ of $H_1(\Sigma_1)$, which we represent as an embedded 
circle in the torus also denoted $\gamma$.

The connected sum theorem for $\widehat{HFK}$ then gives:

\begin{prop}\label{HFKprop} The knot Floer homology of
$\tilde{K}\subset Y = M_n \#^{2g-2} (S^1\times S^2)$ is given by
\begin{equation}
\widehat{HFK}(Y,\tilde{K}) = \Lambda^*H^1(\Sigma_{g})\oplus \left[ 
\Lambda^*H^1(\Sigma_{g-1})\otimes H^*(\coprod_n S^1)\right].
\label{septwistHFK}
\end{equation}
\label{septwistHFKprop}
The only nontrivial differential in the spectral sequence for
$\widehat{HF}(Y)$ is given by contraction with a generator $\gamma\in
H_1(\Sigma_g)$ in the first summand above.
\end{prop}

Indeed, it is shown in \cite{OSknot} that $\widehat{HFK}(B(0,0))\cong
\Lambda^*H^1(\Sigma_1)$ with centered grading. Formula
\eqref{septwistHFK} follows from this and the K\"unneth formula for
$\widehat{HFK}$ under connected sum \cite{OSknot}.

\subsection{Verification of Hypotheses}

Knowledge of the knot Floer homology $\widehat{HFK}(Y,K)$ can often
lead to understanding of the full Heegaard Floer groups $HF^+$ for $Y$
and the surgered manifold $Y_n(K)$. Indeed, there is a spectral
sequence for $HF^{\infty}(Y)$ associated to the filtration of
$CFK^{\infty}(Y,K)$ given by $[\x,i,j] \mapsto i+j$, whose $E_1$-term
is $\widehat{HFK}(Y,K)\otimes\zee[U,U^{-1}]$. The $d_1$ differential
is a sum of homomorphisms that map the group at position $(i,j)$ in
$CFK^\infty$ to those at positions $(i-1,j)$ and $(i, j-1)$: the
``vertical'' and ``horizontal'' components of $d_1$. These components
can in turn be determined from the spectral sequence for
$\widehat{HF}(Y)$ coming from $CFK^{0,*}$: indeed, the vertical
component is precisely (after a translation by a power of $U$) the
first differential in the latter sequence. On the other hand the
complex $CFK^{*,0}$ can also be identified with a filtered version of
$\widehat{CF}$, so the horizontal component of $d_1$ is also
determined by the differential in the spectral sequence calculating
$\widehat{HF}$ from $\widehat{HFK}$. We put these ideas to work in
understanding the Floer homologies of the particular $Y$ and $Y_n$
relevant to our situation.

We begin by determining the differentials in the spectral sequence for
$CF^{\infty}$, in the case of $(Y,\tilde{K})$ as in Proposition
\ref{septwistHFKprop}. As noted above, the $E_1$ term is given by
$\widehat{HFK}(Y,\tilde{K})\otimes\zee[U,U^{-1}]$. Explicitly, this is
\begin{equation}\label{E1term}
E_1 = \left( \Lambda^*H^1(\Sigma_g) \otimes
\zee[U,U^{-1}]\right) \oplus \left( \Lambda^*H^1(\Sigma_{g-1})\otimes
H^*(\coprod_n S^1)\otimes\zee[U,U^{-1}]\right).
\end{equation}
The $d_1$ differential is nontrivial only on the first
summand, where its action is described as follows. Decompose $\Sigma_g
= \Sigma_1\#\Sigma_{g-1}$ where the generator $\gamma\in
H_1(\Sigma_1)$ of Proposition \ref{septwistHFKprop} is contained in
the first factor. Let
\begin{eqnarray*}
E_+ &=& (\Lambda^0H^1(\Sigma_1) \oplus \Lambda^2 H^1(\Sigma_1))\otimes
\Lambda^*H^1(\Sigma_{g-1}) \otimes \zee[U,U^{-1}] \\
E_- &=& \Lambda^1H^1(\Sigma_1)\otimes\Lambda^*H^1(\Sigma_{g-1})
\otimes \zee[U,U^{-1}],
\end{eqnarray*}
so that $E_+\oplus E_-$ is isomorphic to the first summand of the
$E_1$ term above. Then one can check just as in \cite{OSknot} or
\cite{us} that the $d_1$ differential is trivial on $E_+$ while on
$E_-$ it is given by
\begin{equation}\label{d1diff}
d_1(\omega\otimes U^j) =
\iota_\gamma\omega\otimes U^j + PD(\gamma)\wedge\omega\otimes U^{j+1}.
\end{equation}
It is a straightforward exercise to check that the
homology of this differential is given by
\begin{eqnarray*}
H(E_+\oplus E_-, d_1) &=& \Lambda^0H^1(\Sigma_1)\otimes
\Lambda^*(\Sigma_{g-1})\otimes \zee[U,U^{-1}] \\
&&\quad \oplus PD(\gamma)\wedge \Lambda^*H^1(\Sigma_{g-1})
\otimes\zee[U,U^{-1}]\\
&\cong& \Lambda^*H^1(\Sigma_{g-1})\otimes H^*(S^1)\otimes\zee[U,U^{-1}].
\end{eqnarray*}

Therefore the second term in our spectral sequence appears as
\begin{eqnarray*}
E_2 &=& \left(\Lambda^*H^1(\Sigma_{g-1})\otimes
H^*(S^1)\otimes\zee[U,U^{-1}]\right) \\&&\oplus \left(
\Lambda^*H^1(\Sigma_{g-1})\otimes H^*(\coprod_n
S^1)\otimes\zee[U,U^{-1}]\right)\\
&=& \Lambda^*H^1(\Sigma_{g-1})\otimes H^*(\coprod_{n+1} S^1)\otimes 
\zee[U,U^{-1}].
\end{eqnarray*}
Observe that there is an isomorphism of $\zee[U]$-modules
$H^*(S^1\amalg S^1)\otimes\zee[U,U^{-1}] \cong
\Lambda^*H^1(\Sigma_1)\otimes \zee[U,U^{-1}]$. Therefore we can write the
above as
\begin{equation}\label{E2term}
E_2 = \left( \Lambda^*H^1(\Sigma_g)
\otimes\zee[U,U^{-1}]\right) \oplus \left(
\Lambda^*H^1(\Sigma_{g-1})\otimes H^*(\coprod_{n-1}
S^1)\otimes\zee[U,U^{-1}]\right),
\end{equation}
where we have written $\Lambda^*H^1(\Sigma_{g-1})\otimes
\Lambda^*H^1(\Sigma_1) = \Lambda^*H^1(\Sigma_g)$.

We must now determine subsequent differentials in the spectral
sequence, if any. Note first that if $n= 1$ then the above reduces
to $\Lambda^*H^1(\Sigma_g)\otimes \zee[U,U^{-1}]\cong HF^{\infty}(Y)$.
Indeed, this identification follows using the connected sum theorem
for $HF^{\infty}$ (recall that $Y = M_n\#(2g-2)S^1\times S^2$) and the
fact that since $b_1(M_n) = 2$, the Floer homology $HF^\infty(M_n)$ is
``standard'' (see \cite{OS2}, \cite{OS4}). Thus when $n = 1$ there are
no subsequent differentials, and in fact this case has already been
understood in \cite{us}. 

In general, the $d_2$ differential is a sum of three terms, mapping
$C\{i, j\}$ into $C\{i-2, j\}\oplus C\{i-1,j-1\}\oplus C\{i,j-2\}$. In
our case, however, the vertical and horizontal components must be
trivial because those (just as in the case of $d_1$) correspond to
differentials in the second term of the spectral sequence for
$\widehat{HF}$, which we have seen collapses at the second term.
Therefore $d_2$ must be given by a map $d_2: E_2\{i,j\}\to
E_2\{i-1,j-1\}$.

\begin{lemma}\label{sslemma} In terms of the expression \eqref{E2term},
$d_2$ is trivial on the first factor and acts on the second factor by
an isomorphism
\[
\Lambda^*H^1(\Sigma_{g-1})\otimes H^0(\coprod_{n-1} S^1) \otimes U^k
\longto \Lambda^*(\Sigma_{g-1})\otimes H^1(\coprod_{n-1} S^1)\otimes
U^{k+1}.
\]
All subsequent differentials in the spectral sequence are trivial.
\end{lemma}

\begin{proof} We can see that further differentials must be trivial by
examining the gradings. In \eqref{E2term}, the exterior algebras are
equipped with the centered grading, the grading on
$H^*(S^1)$ is the natural homological grading, and $U$ is considered
to have grading $-2$. The filtration (the ``$j$-coordinate'') is equal 
to the grading on the exterior algebras and on
$\zee[U,U^{-1}]$, but $H^*(\coprod_{n-1}S^1)$ is considered to lie in
filtration level 0. (All of these observations can be deduced from the
remarks after Propositions \ref{MnKprop} and \ref{HFKprop}.) Finally, the
``$i$-coordinate'' is recovered by recalling that in the expression
$E_1 = \widehat{HFK}\otimes\zee[U,U^{-1}]$, the subgroup
$\widehat{HFK}\otimes 1$ lies in the column $i = 0$. It is now
straightforward to see that if $(i+j) - (i' + j') > 2$ then there are
no elements $a\in C\{i,j\}$ and $b\in C\{i',j'\}$ whose degrees differ
by 1, so differentials beyond $d_2$ vanish for dimensional reasons.

It follows that the homology of $d_2$ must yield $HF^{\infty}(Y)\cong
\Lambda^*H^1(\Sigma_g)\otimes \zee[U,U^{-1}]$; furthermore the only
factors in $E_2\{i,j\}$ and $E_2\{i-1,j-1\}$ that can be connected by
this differential (i.e., factors whose degrees differ by 1) are those
that are indicated in the statment. Dimensional considerations ensure
that the differential must be an isomorphism between those factors.
\end{proof}

With this understanding of the differentials in the spectral sequence
it is a straightforward matter to determine $HF^+(Y)$ and $HF^+(Y_n)$
to a degree sufficient to verify the hypotheses of Theorem
\ref{0surgthm}. Indeed, hypotheses 1 and 2 of that theorem follow
from:

\begin{lemma} With $Y = M_n\#(2g-2)S^1\times S^2$ as above, we have
an identification of $\zee[U]$-modules
\[
HF^+(Y, \s) = \left(\Lambda^*(\Sigma_{g-1})\otimes \zee[U^{-1}]\right)
\oplus (n-1)\Lambda^*H^1(\Sigma_{g-1}),
\]
where $\s$ denotes the torsion \spinc structure on $Y$, and
$(n-1)\Lambda^*H^1(\Sigma_{g-1})$ denotes the direct sum of $n-1$
copies of the exterior algebra. In particular, $HF_{red}(Y,\s) =
(n-1)\Lambda^*H^1(\Sigma_{g-1})$ as graded groups (where the grading
on the exterior algebra is centered as before), and the filtration on
$HF_{red}(Y,\s)$ is equal to the grading.
\end{lemma}

Note that one can prove this (except for the information about the
filtration) without using the spectral sequence by appealing to the
connected sum theorem for Heegaard Floer homology.

\begin{proof} By restriction, the filtration ${\cal F}:
[\x,i,j]\mapsto i+j$ used to produce the spectral sequence for
$HF^\infty$ also gives a filtration on $CF^+$ and thereby a spectral
sequence for $HF^+$ whose differentials are just the restrictions of
the originals to this quotient complex. In particular the $E_1$ term
appears as $\widehat{HFK}(Y,\tilde{K})\otimes\zee[U^{-1}]$, with
differential given by \eqref{d1diff}. One checks that no new cycles
are created in $E_1$ by passing to the quotient complex $C\{i\geq
0\}$, so that the $E_2$ term here looks just like \eqref{E2term} with
$\zee[U,U^{-1}]$ replaced by $\zee[U^{-1}]$. 

The second differential takes the same form as previously, but in this
case there are additional cycles in $E_2$: since $d_2$ maps
$E_2\{i,j\}$ into $E_2\{i-1,j-1\}$, it sends those elements lying in the
group $\Lambda^*H^1(\Sigma_{g-1})\otimes H^0(\coprod_{n-1}S^1) \otimes
U^0$ (supported in the column $i = 0$) to 0, where it did not do so in the
spectral sequence for $HF^\infty$. Hence this group is precisely
$HF_{red}(Y,\s)$, and can be written $(n-1)\Lambda^*H^1(\Sigma_{g-1})$
as a graded group. The lemma follows immediately, keeping in mind the
structure of the filtration as described in the proof of Lemma
\ref{sslemma}.
\end{proof}

For the remaining hypotheses in Theorem \ref{0surgthm} we must
understand $HF^+(Y_n(\tilde{K}), k)$.
Recall the isomorphism ${}^b\Psi^+: CF^+(Y_n)\cong C\{i\geq 0 \mbox{ or } 
j\geq k\}$ of Theorem \ref{hello}. Using ${}^b\Psi^+$, the filtration
$\cal F$ on $CF^\infty(Y)$ restricts to a filtration on $CF^+(Y_n)$,
and thereby gives a spectral sequence for $HF^+(Y_n)$ It is a simple
matter to see, by examining the domains and ranges of the
differentials, that $HF_{red}(Y_n)$ must be supported along the
right-angled strip $\mbox{max}\{i, j-k\} = 0$. (Indeed, just as in the
case of $Y$ the reduced homology is a result of the ``additional''
cycles for $d_2$ that arise from passing to the quotient complex
$C\{i\geq 0\mbox{ or }j\geq 0\}$, which can only lie in the indicated
region of the $(i,j)$ plane.) Since the $j$-coordinate measures the
filtration, hypothesis 4 follows immediately.

Finally for hypothesis 3, note that the proof of Theorem
\ref{0surgthm} shows that we may replace $F$ by $\pi_*$ since the
corresponding kernels and cokernels are isomorphic once the other
hypotheses hold (c.f. Lemmas \ref{rinvlemma} and \ref{rinvlemma2}). But
it is clear that $\ker \pi_*$ is equal to that portion of
$H_*(C\{i\geq 0 \mbox{ or }j\geq k\})$ lying in the region $i<0$,
while its image is $\im(U^r) \oplus HF_{red}^{\leq k}(Y)\subset
HF^+(Y)$. In particular
\[
\coker (\pi_*) = HF_{red}^{>k}(Y) \cong \Lambda^{>g-1+k}H^1(\Sigma_{g-1})
\]
is a free $\zee$-module.

\subsection{Calculation}\label{proofsec}

We turn our attention to determining the Heegaard Floer homology groups
$HF^+(M(t_\sigma^n), \s)$ where $\s$ is any nontorsion \spinc
structure, $\sigma\subset \Sigma_g$ is a separating curve such that
$\Sigma_g\setminus\sigma$ consists of components of genus $1$ and
$g-1$, $t_\sigma$ denotes the right-handed Dehn twist about $\sigma$,
and $M(t_\sigma^n)$ is the mapping torus of the diffeomorphism
$t_\sigma^n$ for any $n\neq 0$. We focus first on the case of $n>0$.

The result of the preceding section is that the desired Floer homology
can be determined from the knot complex for $\tilde{K}\subset Y$ as
$H_*(C\{i<0\mbox{ and } j\geq k\})$, where $\tilde{K} =
K\#(g-1)B(0,0)$ as before. We assume here that $k>0$. We make use of the same
spectral sequence as in that section to calculate this homology; note
that the differentials have already been determined. Now, the $E_1$
term of this spectral sequence is just the portion of
$\widehat{HFK}(Y,\tilde{K})\otimes \zee[U,U^{-1}]$ that lies in the
relevant part of the $(i,j)$ plane. Following \cite{OSknot} and
\cite{us}, we introduce the notation
\begin{equation} \label{Xgddef}
X(g,d) = \bigoplus ^d _{i = 0} \Lambda ^{2g-i} H^1(\Sigma _g) \otimes
_\mathbb{Z} \mathbb{Z}[U] /U^{d+1-i}.
\end{equation}
Thus $X(g,d) \cong H^*(\sym^d\Sigma_g)$ as $\zee[U]$-modules (see
\cite{macdonald}). Then it is easy to see that in the spectral
sequence for $HF^+(Y_0(K), \s_k)$,
\[
E_1 = X(g,d) \oplus \left[ X(g-1,d-1)\otimes H^*(\coprod_n 
S^1)\right],
\]
where $d = g-1-k$. The $d_1$ differential is nontrivial only on the first factor (see
the discussion after equation \eqref{E1term}), where it acts as in
equation \eqref{d1diff}. The complex $(X(g,d), d_1)$ was considered in
\cite{us}, where it was shown that its homology is
\[
H_*(X(g,d), d_1) = (X(g-1,d-1)\otimes H^*(S^1))\oplus
\Lambda^{2g-2-d}H^1(\Sigma_{g-1})_{(g-d)}.
\]
Hence, 
\begin{eqnarray}
E_2 &=& [X(g-1,d-1)\otimes H^*(S^1)\oplus
\Lambda^{2g-2-d}H^1(\Sigma_{g-1})_{(g-d)}] \nonumber\\
&&\oplus \left[ X(g-1,d-1)\otimes H^*(\coprod_n S^1)\right]\nonumber\\
&=& \left[ X(g-1,d-1)\otimes H^*(S^1\sqcup S^1)\oplus
\Lambda^{2g-2-d}H^1(\Sigma_{g-1})_{(g-d)}\right] \label{E2term2} \\
&&\oplus \left[ X(g-1,d-1)\otimes H^*(\coprod_{n-1} S^1)\right].\nonumber
\end{eqnarray}
The $d_2$ differential acts as was determined in Lemma \ref{sslemma},
only on the last summand above (and only when $n>1$). Thus in ``most''
positions $(i,j)$ the last factor is killed in homology, with the
exception of those $(i,j)$ with $i = 0$ or $j = k$: since the
differential maps $(i,j)$ to $(i-1,j-1)$ there are additional cycles
when $j = k$ and fewer boundaries when $i = 0$. Specifically, the
homology of \eqref{E2term2} is given by the first term in brackets
plus the contributions:
\begin{eqnarray}
(i =0 ) &\quad& \bigoplus_{p = 1}^{d} \Lambda^{2g-2-d + p}
H^1(\Sigma_{g-1})\otimes H^1(\coprod_{n-1} S^1)\label{vertcont}\\
(j = k) &\quad& \bigoplus_{p = 1}^{d}
\Lambda^{2g-2-d+ p}H^1(\Sigma_{g-1}) \otimes H^0(\coprod_{n-1} S^1)\otimes
U^{p-1}.\label{horizcont}
\end{eqnarray}
where $d = g-1-k$.

\vspace{1em}
\noindent
{\it Proof of Theorem \ref{mainthm}}.
The results of the preceding sections show that
for given $k$ with $0<k\leq g-1$
\[
HF^+(M(t_\sigma^n), \s) \cong H_*(C\{i<0\mbox{ and } j\geq k\})
\]
where the right-hand side refers to the homology of the indicated
quotient complex of $CFK^\infty(Y, \tilde{K})$ (in the torsion \spinc
structure on $Y = M_n\# (2g-2)S^1\times S^2$). The latter has been
shown to be isomorphic as a relatively graded $\zee$-module to the sum of
the first bracketed term in \eqref{E2term2} and the two expressions
\eqref{vertcont} and \eqref{horizcont}.

To understand those two expressions, recall that the summand
$\Lambda^{2g-q}H^1(\Sigma_g)\otimes U^j$ of $X(g,d)$ is supported in
grading $g-q-2j$. From this it follows (recall that the factors in
\eqref{vertcont} and \eqref{horizcont} arise from $X(g-1,d-1)$ as in
the second term of \eqref{E2term2}) that
\begin{eqnarray*}
\mbox{\eqref{vertcont} $\oplus$ \eqref{horizcont}} &\cong& (n-1)
\bigoplus_{p = 1}^d \Lambda^{2g-2-d+p}_{(g-d+p)}H^1(\Sigma_{g-1}) \oplus
\Lambda^{2g-2-d+p}_{(g-d-p+1)}H^1(\Sigma_{g-1}) \\
&\cong& \bigoplus_{p=1}^d
\Lambda^{2g-2-d+p}_{(g-d-p+1)}H^1(\Sigma_{g-1}) \otimes
H_*(\coprod_{n-1} S^{2p-1}).
\end{eqnarray*}
Adding this to the first line of \eqref{E2term2} gives \eqref{answer},
when $k$ is positive. For negative $k$, the result follows from the
conjugation invariance of $HF^+$ (see \cite{OS1}).

\subsection{The Case of Left-Handed Twists}\label{lefttwist}

The procedure for calculation of $HF^+(M(t_\sigma^{-n}))$, $n>0$, 
is very similar to the positive-twist case. We outline here the main 
differences. 

The surgery diagram for $M(t_\sigma^{-n})$ is identical to that for 
$M(t_\sigma^n)$ with the exception that $-1$-surgery curves 
corresponding to the Dehn twists are replaced by $+1$-curves (c.f. 
Figure \ref{sepknot2}). Thus, $M(t_\sigma^{-n})$ is obtained by 
0-surgery on the knot $\tilde{K} = \#(g-1)B(0,0)\# K$ as before, 
where $K\subset M_{-n}$ is as in Figure \ref{sepknot}(a) with 
signs changed, or as in Figure \ref{sepknot}(b) with positive twists 
replaced by negative ones. Since $M_{-n}$ is just $M_n$ with the 
orientation reversed, we have (using the behavior of $\widehat{HF}$ 
under orientation reversal) $\widehat{HF}(M_{-n}; \s_0) = 
\zee^{n+1}_{(0)} \oplus \zee^{n+1}_{(-1)}$.

To obtain the knot Floer homology $\widehat{HFK}(M_{-n}, K)$, 
we proceed as before: the result corresponding to Lemma \ref{M1Klemma} 
is
\[
\widehat{HFK}(M_{-1},K; j) = \left\{\begin{array}{ll} \zee_{(1)} & j = 
1\\ \zee_{(-1)} \oplus \zee^3_{(0)} & j = 0 \\ \zee_{(-1)} & j = -1 \\ 
0 & \mbox{otherwise.}\end{array}\right.
\]
This is proved in just the same way as Lemma \ref{M1Klemma} with some 
surgery coefficients having opposite sign (see \cite{us}).

In general, the analogue of Proposition \ref{MnKprop} is:
\[
\widehat{HFK}(M_{-n}, K; j) = \left\{\begin{array}{ll}\zee_{(1)} & j = 
1\\ \zee^n_{(-1)} \oplus \zee^{n+2}_{(0)} & j = 0 \\ \zee_{(-1)} & 
j = -1 \\ 0 & \mbox{otherwise,}\end{array} \right.
\]
and is proved in an analogous manner. Put another way,
\[
\widehat{HFK}(M_{-n}, K; j) \cong \Lambda^*H^1(\Sigma_1)\oplus 
H^*(\coprod_{n} S^1)[-1],
\]
where ``$[-1]$'' indicates that the grading on the corresponding 
factor has been shifted down by 1 (the grading on the exterior 
algebra, as usual, is centered). Note that the differential $d_1$ in 
the spectral sequence calculating $\widehat{HF}(M_{-n})$ from the 
above consists of an injection $\zee_{(1)}\to \zee^{n+2}_{(0)}$, 
which we may identify with contraction by a generator in 
$H_1(\Sigma_1)$. 

It follows just as in the previous work that 
\[
\widehat{HFK}(Y,K) \cong \Lambda^*H^1(\Sigma_g) \oplus \left( 
\Lambda^*H^1(\Sigma_{g-1})\otimes H^*(\coprod_n S^1)[-1]\right),
\]
and in the spectral sequence calculating $HF^+(M(t_\sigma^{-n}))$,
\[
E_1 = X(g,d) \oplus \left( X(g-1,d-1)\otimes H^*(\coprod_n 
S^1)[-1]\right).
\]
The $d_1$ differential in this case again acts only on the first 
factor above, and is given by equation \eqref{d1diff}. However, under 
the splitting $E_+\oplus E_-$ of the first factor (c.f. the remarks 
after equation \eqref{E1term}), here $d_1$ is trivial on $E_-$ and 
acts nontrivially only on $E_+$. The homology of $X(g,d)$ with 
respect to this differential was also calculated in \cite{us}, and 
that result gives
\begin{eqnarray*}
E_2 &=& X(g-1,d-1)\otimes H^*(S^1)[-1] \oplus 
\Lambda^{2g-2-d}H^1(\Sigma_{g-1})_{(g-d)} \\
&& \oplus \left( X(g-1,d-1)\otimes H^*(\coprod_n S^1)[-1]\right)\\
&\cong& X(g-1,d-1)\otimes H^*(S^1\sqcup S^1)[-1] \oplus 
\Lambda^{2g-2-d}H^1(\Sigma_{g-1})_{(g-d)} \\
&& \oplus \left( X(g-1,d-1)\otimes H^*(\coprod_{n-1} S^1)[-1]\right).
\end{eqnarray*}

An argument as in the positive case identifies the $d_2$ differential 
as being nontrivial only on the second factor above, mapping
\[
d_2: \Lambda^kH^1(\Sigma_{g-1})\otimes 
H^0(\coprod_{n-1}S^1)[-1]\otimes U^j \longrightarrow \Lambda^kH^1(\Sigma_{g-1})\otimes 
H^1(\coprod_{n-1}S^1)[-1]\otimes U^{j+1}
\]
isomorphically. The resulting homology is given again by the sum of 
\eqref{vertcont} and \eqref{horizcont}, this time with a grading 
shift of $-1$. This proves Theorem \ref{mainthm} in the case $n<0$.


\end{document}